\documentstyle[12pt]{article}

  \newcommand{\const}{\rm const}
  \newcommand{\Var}{\rm Var}

\begin{document}

   \begin{center}

 {\bf  Moment and exponential tail estimations for norms of random variables} \par

 \vspace{3mm}

{\bf  and random operators }\par

\vspace{3mm}

 {\bf in mixed (anisotropic) Lebesgue - Riesz spaces. } \par

\vspace{5mm}

{\bf M.R.Formica,  E.Ostrovsky and L.Sirota}.\par

\end{center}

 Universit\`{a} degli Studi di Napoli Parthenope, via Generale Parisi 13, Palazzo Pacanowsky, 80132,
Napoli, Italy. \\

e-mail: mara.formica@uniparthenope.it \\

\vspace{3mm}

Department of Mathematics and Statistics, Bar-Ilan University, \\
59200, Ramat Gan, Israel. \\

e-mail: eugostrovsky@list.ru\\

\vspace{3mm}

\ Department of Mathematics and Statistics, Bar-Ilan University, \\
59200, Ramat Gan, Israel. \\

e-mail: sirota3@bezeqint.net \\

\vspace{5mm}

\begin{center}

 \ {\bf Abstract.} \par

\end{center}

\vspace{5mm}

 \hspace{3mm} We study the random variables  (r.v.)  with values in the so - called mixed (anisotropic) Lebesgue - Riesz spaces:
 formulate the sufficient conditions for belonging  of the r.v. to these spaces, estimate the tail of norms distribution, especially
 deduce the exponential decreasing tails of them, etc.\par
  \ We obtain as a consequence the estimations of the norms of {\it random} integral operators acting between these spaces. \par

\vspace{5mm}

\begin{center}

{\sc Key words and phrases.}

\end{center}

\vspace{5mm}

 \ Random variables and vectors (r.v.), distribution and tails of distributions, Tchebychev - Markov inequality, moment, kernel, H\"older's
inequality,  Lebesgue - Riesz and Grand Lebesgue Spaces (GLS), norms, Young - Fenchel transform, Young inequality, exponential decreasing function,
random and ordinary  linear integral  operators, measure, mixed (anisotropic) spaces, permutation inequality, factorable functions.

\vspace{5mm}

\section{Statement of problem.  Notations. Previous works.}

\vspace{5mm}

 \hspace{3mm}  We recall here the definition  of the so-called anisotropic Lebesgue (Lebesgue-
Riesz) spaces, which appeared in the famous article of Benedek A. and Panzone R.  \cite{Benedek Panzone}.
 More detail information about this spaces with described applications see in the classical
books of Besov O.V., Ilin V.P., Nikolskii S.M. \cite{Besov Ilin Nikolskii}, chapter 1,2; Leoni G. \cite{Leoni}, chapter 11.\par

\vspace{3mm}

 \ Let $ \ (\Omega = \{\omega\}, {\cal B}, {\bf P}) \ $ be certain probability space with expectation $ \ {\bf E } \ $
 and variation $ \ {\bf \Var,} \ $ and let \\

$$
 (X_k =\{x_k\}, A_k, \mu_k), \ k = 1,2,\ldots, l, \ l < \infty
$$

be measurable spaces with sigma - finite separable non - trivial measures
$ \ \mu_k. \ $  The separability denotes that the metric space  $ \ A_k \ $ relative the distance

$$
\rho_k(D_1,D_2) = \mu_k(D_1 \setminus D_2) + \mu_k(D_2 \setminus D_1) = \mu_k(D_1 \Delta D_2), \ D_{1,2} \subset A_k
$$
is separable. \par

 \ Denote  $ \  X := \otimes_{k=1}^l X_k,  \ $ so that $ \ x \in X \ \Leftrightarrow \ x = \vec{x} = \{\ x_1,x_2,\ldots, x_l \ \}.  \ $  \par

 \vspace{3mm}

 \ Let also $ \ \vec{p} =  p = (p_1, p_2, ..., p_l) \ $ be $ \ l \ $ dimensional fixed numerical vector such that $ \  1 \le p_j < \infty. \ $
Recall that the anisotropic (mixed) Lebesgue - Riesz space

\begin{equation} \label{anisotr space X}
 L(p) = L_p = L(\vec{p}) = L(\vec{p}, \ \{X_k\}, \ \{\mu_k\})
\end{equation}
consists  by definition on all the
total measurable real valued function $ \  f = f(x_1, x_2, . . . , x_l) = f(\vec{x}): \ $

$$
f: \otimes_{k=1}^l X_k \to R,
$$
having a finite norm  $ \ ||f||_{\vec{p}} \stackrel{def}{=} \ $

\vspace{3mm}

\begin{equation} \label{def mixed norm}
\left( \int_{X_l} \mu_l(dx_l) \left(\int_{X_{l-1}} \mu_{l-1}(dx_{l-1}) \ldots \left(\int_{X_1} \mu_1(dx_1) \ |f(\vec{x})|^{p_1} \right)^{p_2/p_1} \right)^{p_3/p_2} \ldots  \ \right)^{1/p_l}.
\end{equation}

\vspace{3mm}

 \ In particular, for the  one - dimensional numerical valued r.v. $ \  \xi = \xi(\omega) \ $ as well as for the number $ \ p, p \in [1,\infty) \ $ we obtain the
 classical Lebesgue - Riesz $ \ L_p(\Omega) \ $ norm

 $$
 ||\xi||_p = \left( \ {\bf E} |\xi|^p \ \right)^{1/p},
 $$
as well as ones

$$
 ||f||_{p,X_k} = \left[ \int_{X_k} \ |f(x_k)|^p \ \mu_k(dx_k) \ \right]^{1/p}.
$$

 \vspace{3mm}

 \ Note that in general case $ \  ||f||_{p_1,p_2} \ne ||f||_{p_2,p_1}, $  but $ \ ||f||_{p,p} = ||f||_p. \ $ \par

\vspace{3mm}

 \ Observe also that if $ \ f(x_1, x_2) = g_1(x_1) \cdot g_2(x_2), \ $  (condition of factorization), then

\begin{equation} \label{factorization}
||f||_{p_1,p_2} = ||g_1||_{p_1} \cdot ||g_2||_{p_2},
\end{equation}

(formula of factorization).\par

\vspace{3mm}

 \ Note that under conditions separability  of measures $ \ \mu_k $  these spaces are also separable and Banach spaces.\par

 \ These spaces arises in the Theory of Approximation, Functional Analysis, theory
of Partial Differential Equations, theory of Random Processes etc.\par

\vspace{3mm}

 \ Let for example $ \ l = 2;\ $  we agree to rewrite for clarity the expression for $ \ ||f||_{p_1,p_2} \ $
as follows:

$$
||f||_{p_1,p_2} := ||f||_{p_1,X_1; \ p_2,X_2}.
$$

 \ Analogously,

$$
||f||_{p_1,p_2,p_3} = ||f||_{p_1,X_1; \ p_2,X_2; \ p_3,X_3}.
$$

 \vspace{3mm}

 \ Notice that the last expression may be rewritten as follows:

\begin{equation} \label{order norms}
||f||_{p_1,p_2,p_3} =  || \ || \ ||f||_{p_1,X_1} \ ||_{p_2, X_2} \ ||_{p_3,X_3}.
\end{equation}

 \vspace{3mm}

 \ Let us recall also the following important fact: the so - called {\it permutation} inequality, in the terminology of an  original article
\cite{Benedek Panzone};  see also  the monograph \cite{Besov Ilin Nikolskii}, chapter 1, pp. 24 - 26. Indeed, let  $ \ (Z,B, \nu) \ $  be {\it another}
measurable space and $ \ \phi: \ (X,Z) = X \otimes Z \to R \ $ be common measurable function. In what follows as before $ \ X = \otimes_{k=1}^l X_k. \ $ \par
 \  Let also $ \ r = \const \ge \overline{p}, \ $ where

$$
\overline{p} := \max_{k=1}^d p_k.
$$

 \  It is true the following {\it permutation} inequality (in our notations):

\vspace{3mm}

\begin{equation} \label{permutation}
||\phi||_{p,X; \ r,Z}  \le ||\phi||_{r,Z; \ p,X}.
\end{equation}

\vspace{3mm}

 \ In what follows $ \  Z = \Omega, \ \nu = {\bf P}. \ $ \par

\vspace{5mm}

\section{Estimation of the distribution for the norm of random field.}

\vspace{5mm}

 \hspace{3mm} Let $ \  \eta = \eta(x,\omega), \ x \in X \ $ be separable total measurable numerical valued random field (r.f.).
 Define the following random variable

\begin{equation} \label{estim norm rv}
\zeta = \zeta(\omega) \stackrel{def}{=} ||\eta||_{\vec{p}} = ||\eta||_{\vec{p},X}.
\end{equation}

\ We intent to estimate  in this section the Lebesgue - Riesz probabilistic norms of the r.v. $ \ \zeta: \ |\zeta|_{r,\Omega}, \ r \ge 1. \ $

\vspace{5mm}

\ {\bf Theorem 2.1.} \ Denote

$$
\psi[p](r) \stackrel{def}{=} || \ ||\eta||_{r,\Omega} \ ||_{p,X}.
$$

 \ Suppose  $ \  r \ge \overline{p};  \ $  then

\vspace{3mm}

\begin{equation} \label{norm estim}
\left[{ \ \bf E} |\zeta|^r  \ \right]^{1/r}  = ||\zeta||_{r,\Omega} \le   \psi_p(r).
\end{equation}

\vspace{4mm}

 \ {\bf Proof} is simple.  We have

$$
||\zeta||_{r,\Omega} = || \  ||\eta||_{\vec{p, X}}  \ ||_{r,\Omega}.
$$

 \ One can apply the permutation inequality (\ref{permutation}):

$$
|| \  ||\eta||_{\vec{p, X}}  \ ||_{r,\Omega} \le || \  ||\eta||_{r,\Omega} \ ||_{\vec{p},\Omega} = \psi[p](r),
$$
 Q.E.D. \par

 \vspace{4mm}

 \  Suppose in addition that the introduced above function $ \ \psi[\vec{p}](r) \ $ is finite for some non - trivial segment
 $ \ r     \in [a,b),   \  $ where $ \ a = \overline{p}, \  a <  b = \const \le \infty.  \ $  The proposition of theorem (2.1)
 may be reformulated as follows. Introduce as ordinary the so - called Grand Lebesgue Space $ \ G\psi_p \ $ builded on our probability space
 consisting on all the r. v.- s $  \xi \ $ having a finite norm

\begin{equation}  \label{grand lebesgue norm}
||\xi||G\psi[p](\cdot)  \stackrel{def}{=} \sup_{r \in (a,b)} \ \left\{ \ \frac{||\xi||_r}{\psi[p](r)} \ \right\}.
\end{equation}

 \ The general theory of these spaces is represented in many works, see e.g. \cite{Ahmed Fiorenza Formica at all}, \cite{anatriellofiojmaa2015},
\cite{anatrielloformicaricmat2016},  \cite{Buldygin-Mushtary-Ostrovsky-Pushalsky}, \cite{caponeformicagiovanonlanal2013},  \cite{Ermakov etc. 1986},
\cite{Fiorenza2000}, \cite{fiokarazanalanwen2004}, \cite{fioguptajainstudiamath2008}, \cite{Fiorenza-Formica-Gogatishvili-DEA2018},
\cite{fioforgogakoparakoNAtoappear}, \cite{fioformicarakodie2017},\cite{Kozachenko}, \cite{Ostr Sir CLT  mixed},\cite{Ostrovsky 0},
\cite{Ostrovsky 1}, \cite{Ostrovsky HAIT} etc. In particular, these spaces are complete, Banach functional and rearrangement invariant.\par

\vspace{4mm}

\ We have from the proposition  (\ref{norm estim})   of theorem 2.1

\vspace{3mm}

\begin{equation} \label{gpsi norm estim}
||\zeta||G\psi[p](\cdot) \le 1.
\end{equation}

 \vspace{3mm}

  \ Introduce the following function

$$
g_p(u) := \sup_{r \in (a,b) } (u r - r \ln \psi[p](r)) \ -
$$
the (regional) Young - Fenchel transformation for the function $ \ r \to r \ln \psi[p](r), \  $
relative the variable $ \ r, \ $  for the values $ \ r \in (a,b). \ $ \par

\vspace{3mm}

 \ It follows from (\ref{gpsi norm estim}) the following {\it exponential decreasing} (in general case) tail estimation for the
 distribution of the random variable $ \ \zeta: \ $

$$
{\bf P} (|\zeta| > u) \le \exp( \ - g_p(u)  \ ), \ u \ge 1.
$$

 \ The inverse conclusion also holds true under appropriate natural conditions, see  \cite{Ermakov etc. 1986}, \cite{Kozachenko},
\cite{Ostrovsky HAIT}. \par

\vspace{4mm}

 \ {\bf  Example 2.1.}  \ Let $ \ \exists \  C(p) < \infty \ \forall r \ge 1 \ \Rightarrow  \psi[p](r) \le C(p) \ r^{1/m}, \  \ m = \const > 0;  \ $  then
 $ \ \exists \ C_2(m,p) > 0 \ \Rightarrow  \ $

$$
{\bf P} (|\zeta| > u) \le \exp  \left( \ - C_2(m,p) \ u^m  \ \right), \ u \ge 1,
$$
and inverse conclusion is also true. \par

\vspace{4mm}

 \ {\bf  Example 2.2.} \ Let us clarify slightly  the applicability of our theorem yet in the case $ \  d = l = 1. \ $ Namely, let again
 $ \  \eta = \eta(x) = \eta(x,\omega)   \ $  be separable measurable numerical valued random field; put

 $$
 \zeta = \zeta(\omega) = ||\eta||_{p,X}, \ p \in [1,\infty).
 $$
 \ How one can  estimate the tail of distribution of the r.v. $ \ \zeta ? \ $ \par

 \vspace{3mm}

  \  Answer. Let $ \ r \ $ be arbitrary number greatest than $ \ p: \ r \ge p. \ $ We propose by virtue of Theorem 2.1

$$
\left[ \ {\bf E}|\zeta|^r  \ \right]^{1/r} \le  \left\{ \ \int_X  \ \mu(dx) \ \left[ \ {\bf E} |\eta(x)|^r \ \right]^{p/r} \ \right\}^{1/p}.
$$

\vspace{5mm}

\section{Estimation of the norm of random  operators.}

\vspace{5mm}

 \ Let also

\begin{equation} \label{anisotrop space Y}
 L(q) = L_q = L(\vec{q}) = L(\vec{q}, \ \{Y_j\}, \ \{\nu_j\}), \ j = 1,2,\ldots, d; \ d < \infty
\end{equation}
be {\it another}  $ \ "d" \ - \ $ dimensional  mixed (anisotropic) space and  put as before \\
 \ $ \  Y := \otimes_{j=1}^d Y_j.  \ $ \par

\vspace{3mm}

 \ Let $ \  K  = K(x, y, \omega), \ y \in Y, \ x \in X, \ \omega \in \Omega \ $ be certain {\it \ kernel,} i.e. a numerical valued
 total measurable function, on the other words, random kernel. The function $ \ K(\cdot, \ \cdot, \ \cdot) \ $ is named also a random field  (r.f.). \par

 \vspace{4mm}

 \ Introduce the following  important linear {\it random} integral operator  $ \ U  \ $  having the kernel $ \ K: \ $

\begin{equation} \label{oper}
U[g](x)  \stackrel{def}{=} \int_Y K(x, y, \omega) \ g(y) \ \prod_{j=1}^d \nu_j(dy_j), \ y \in Y, \ x \in X,
\end{equation}

\vspace{4mm}

or briefly

\begin{equation} \label{oper briefly}
U[g](x)  \stackrel{def}{=} \int_Y K(x, y, \omega) \ g(y) \ \nu(dy), \ y \in Y, \ x \in X, \ \omega \in \Omega.
\end{equation}

\vspace{4mm}

 \ {\bf  Our target in this section is investigation  of conditions for the finiteness a.e. of this operator acting
 between two mixes Lebesgue spaces and estimate the distribution of its norm. } \par

\vspace{4mm}

 \ There are several publications devoted the theme of the random operators acting between  different Banach spaces: \cite{Arnold},
 \cite{Dorogovtsev}, \cite{Howard}, \cite{Thang1}, \cite{Thang2}, a classical  monograph \cite{Skorohod} and so one. The case of stochastic integral
 operators ic considered in particular in  \cite{Arnold}, \cite{Skorohod}.\par

 \vspace{5mm}

 \ Let us recall the following simple estimation for the norm of linear {\it deterministic} integral kernel  operator  of the form

\begin{equation} \label{determin oper}
f(x) = V[g](x) \stackrel{def}{=} \int_Y V(x,y) \ g(y) \ \nu(dy), \hspace{3mm}  x \in X, \ y \in Y.
\end{equation}

 \ It follows immediately by virtue of H\"older's inequality that

\begin{equation} \label{auxiliary determin}
||f||_{r,X} \le || \ || \ V(\cdot,\cdot) \ ||_{q,Y} \ ||_{r,X} \times ||g||_{p,Y},
\end{equation}
where as before $ \  p  = \vec{p} = \{ p_j  \},   \  q = \vec{q} = \{q_j \}, \ i = 1,2,\ldots,d;    \ $

\begin{equation} \label{pqr}
\frac{1}{p_j} + \frac{1}{q_j} = 1, \ p_j, q_j \in (1,\infty),
\end{equation}
and $ \ r = \vec{r} = \{ \ r_k \ \},  \ k = 1,2,\ldots,l; \ 1 < r_k < \infty; \ $
see e.g. the classical monograph \cite{Kantorowitch Akilov}, chapter XI, section 3.\par

\vspace{4mm}

 \hspace{3mm} {\bf Remark 3.1.} There are many works devoted to the operators norm estimates for {\it concrete} form of  these
linear operators: Fourier,  Laplace, Pseudo - differential, singular, convolution, Bochner, Riesz, fractional, Hardy - Littlewood,
Hausdorff etc. The fundamental investigation of this theory may be found in particular in the famous book of  G.O.Okikiolu  \cite{Okikiolu};
see also   the brief review in an article \cite{Ostrovsky 1}.\par

\vspace{4mm}

 \ Let us return to the random integral operator $ \ U, \ $  described above  in (\ref{oper}), (\ref{oper briefly}). We have by virtue of
 (\ref{auxiliary determin}) under at the same notations

\begin{equation} \label{U without omega}
||U[g]||_{r,X} \le || \ || K ||_{q,Y} \ ||_{r,X} \ || \times || g||_{p,Y}.
\end{equation}

 \ Further, let $ \ s \ $ be some number such that

\begin{equation} \label{restric s}
s \ge  \max\{ \ \max_k r_k, \ \max_j q_j \  \}.
\end{equation}

 \ We have taking the norm $ \ ||\cdot||_{s,\Omega} \ $ from both the sides of the inequality (\ref{U without omega}) taking into account
 the permutation inequality

\begin{equation} \label{U with omega}
|| \ ||U[g]||_{r,X} \ ||_{s,\Omega} \le || \ || \ || K ||_{s,\Omega}  \  ||_{q,Y} \ ||_{r,X} \ \times || g||_{p,Y}.
\end{equation}

\vspace{4mm}

 \ To summarize. Denote for these values of the parameters

\begin{equation} \label{theta fun}
\theta(s) = \theta_{p,r}(s) := || \ || \ || K ||_{s,\Omega}  \  ||_{q,Y} \ ||_{r,X},
\end{equation}

\begin{equation} \label{AB}
A := \max\{ \ \max_k r_k, \ \max_j q_j \  \}.
\end{equation}

\vspace{4mm}

 \ {\bf Theorem 3.1.} Suppose that for some values $ \ a \ge A, \ b \in (a,\infty] \ \Rightarrow \ \theta(s) < \infty. \  $  Then

\begin{equation} \label{main oper}
|| \ ||U[g] \ ||_{r,X} \ ||G\theta_{a,b}  \le 1,
\end{equation}

with correspondent tail estimation. \par

\vspace{4mm}

\ {\bf Example 3.1.} Assume that the random kernel $ \ K(x,y,\omega) \ $  allows a {\it factorization:}

$$
| K(x,y,\omega) | \le K_0(x,y) \ \tau(\omega),
$$
where the non - negative (measurable)  functions $ \ K_0(\cdot, \cdot), \ \tau  \ $ are such that

$$
 h(q,r) := || \ || K_0||_{q,Y} \ ||_{r,X} < \infty
$$
and

$$
\exists (a,b) = \const, \ 1 \le a < b \le \infty,  \ \forall s \in (a,b) \ \Rightarrow  \ \rho(s) := ||\tau||_{s,\Omega} < \infty.
$$

 \ Then

\begin{equation} \label{factoriz example}
 \theta_{p,r}(s) \le h(q,r) \cdot \rho(s), \ s \in (a,b).
\end{equation}

\vspace{5mm}

\section{Concluding remarks.}

 \hspace{3mm} {\bf A.} The essential unimprovability of obtained estimations may be illustrated, for instance, in the
 one - dimensional case $ \ l  = d = 1, \ $ see e.g. \cite{Okikiolu}, chapters 3,4;  \cite{Ostrovsky 1}. Namely,
 define the factorable r.f. $ \ \eta(x,\omega) = \tau(\omega) \cdot h(x), \ \omega \in \Omega, \ x \in X, \ $
where for definiteness

$$
||\tau||_{r,\Omega} = 1 = ||h||_{p,X}.
$$
 \ Then both the sides of inequality (\ref{norm estim}) are equal to 1, by virtue of (\ref{factorization}):

$$
 || \   || h \cdot \tau ||_{p,X} \ ||_{r,\omega} = ||h(\cdot)||_{p,X} \times ||\tau||_{r,\Omega} = 1 =
$$

$$
 || \   || h \cdot \tau ||_{r,\Omega} \ ||_{p,X},
$$
still without  the restriction $ \ r \ge p. \ $ \par

 \vspace{4mm}

 \ {\bf B.} Let us ground the unimprovability  of the integral operators norm estimations (\ref{determin oper}),
(\ref{auxiliary determin}). Indeed, assume  as above  $ \  d = l = 1;  \ \nu(Y) = 1,$  and set $ \ g(y) = 1; \ $
and let the kernel $ \ V(\cdot) \ $ be degenerate: $ \ V(x,y) = v(x). \ $ Then both the hand sides of (\ref{determin oper})
are equal:

$$
f(x) = v(x), \ ||f||_{r,X} =||v||_{r,X} = || \ V(\cdot,\cdot) \ ||_{q,Y} \ ||_{r,X} = || \ V(\cdot,\cdot) \ ||_{q,Y} \ ||_{r,X} \times ||g||_{p,Y}.
$$

\vspace{4mm}

\ {\bf C.} The  unimprovability of the assertion of theorem 3.1. about the random linear operators follows formally from ones for 
deterministic operator. One can use also the slightly modified  previous  example, choosing for instance 

$$
K(x,y,\omega) := V(x) \cdot \xi(\omega),
$$
where the non - zero r.v. $ \ \xi \ $ belongs to the Lebesgue - Riesz space $ \ L_r(\Omega), r > 1, \ $ 
and as before $ \ \nu(Y)  = 1, \ g(y) = 1. \ $      

\vspace{4mm}

 \ {\bf D.} The  unimprovability of both the propositions of theorem 2.1 and 3.1 in the multidimensional 
 case may be shown by means of consideration of the {\it factorable functions,} i.e. when

 $$
 g(y) = \prod_{j=1}^d g_j(y_j), \ f(x) = \prod_{k=1}^d f_k(x_k), \ V(x,y) = \prod_{k=1}^d V_k(x_k)
 $$
etc.\par

\vspace{3mm}

  \ {\bf E.} \ It is interest in our opinion to investigate  analogously the case of non - linear random operators,
  as well as consider a possibility when $ \ l = \infty \ $ or $ \ d = \infty. \ $ \par

\vspace{6mm}

\vspace{0.5cm} \emph{Acknowledgement.} {\footnotesize The first
author has been partially supported by the Gruppo Nazionale per
l'Analisi Matematica, la Probabilit\`a e le loro Applicazioni
(GNAMPA) of the Istituto Nazionale di Alta Matematica (INdAM) and by
Universit\`a degli Studi di Napoli Parthenope through the project
\lq\lq sostegno alla Ricerca individuale\rq\rq .\par

\end{document}